\documentclass[12pt,a4paper]{article}
\usepackage{epsfig,amstext,amssymb,amscd}

\begin{document}

\title{Discrete convexity and unimodularity. I}
\author{ Vladimir I.
Danilov and Gleb A. Koshevoy }

\date{}

\maketitle

\section{Introduction}

In this paper we develop a theory of convexity for the lattice of
integer points $\mathbb Z^n$, which we call theory of {\em
discrete convexity}.

What subsets $X\subset {\mathbb Z}^n$ could be called  "convex"?
One property seems indisputable: $X$ should coincide with the set
of all integer points of its convex hull $\mbox{co}(X)$. We call
such sets {\em pseudo-convex}. The resulting class ${\mathcal PC}$
of all pseudo-convex sets is stable under intersection but not
under summation. In other words, the sum $X+Y$ of two
pseudo-convex sets $X$ and $Y$ needs not be pseudo-convex. We
should consider subclasses of $\mathcal {PC}$ in order to obtain
stability under summation.

As we show stability under summation is closely related to another
question: when the intersection of two integer polytopes is an
integer polytope? Beginning from the paper \cite{Ed}, it is known
that the class of generalized polymatroids has this property. Let
us define a $PM$-set in $\mathbb Z^n$ as the set of integer points
of some (integer) g-polymatroid. Then the class of all $PM$-sets
is a class of discrete convexity (DC-class). Specifically, the sum
of {\em PM}-sets is a {\em PM}-set, and non-intersecting {\em
PM}-sets can be separated by some linear functional.

On this way at least two questions arise:

    1) Can we extend the class of g-polymatroids without losing
in the process the nice properties which precisely made us
consider it at the very beginning?

    2) Do other classes exist which exhibit similar properties?
If so, how are they to be constructed or described?

Answers on these questions (`No' on the first one and `Yes' on the
second) rest on a relation of the discrete convexity with
unimodular systems. The latter are nothing but invariant versions
of totally unimodular matrices (we discuss their properties in
Section 5). Every unimodular system $\mathcal R$ defines a class
$\mathcal{P}t(\mathcal{R},\mathbb{Z})$ of integer polytopes which possesses two
properties:

a) it is stable under summation;

b) the intersection of any two polytopes from
$\mathcal{P}t(\mathcal{R},\mathbb{Z})$ is an integer polytope.

The class $\mathcal{P}t(\mathcal{R},\mathbb{Z})$ consists of those integer polytopes
all edges of which are parallel to some elements of $\mathcal{R}$. Moreover, any
ample class of integer polytopes with the properties a) and b) has such a form.

For example, the class of g-polymatroids corresponds to the
unimodular system $\mathbf{A}_n$ in $\mathbf{Z}^n$ which consists of
vectors $\pm e_i$ and $e_i-ej$, $i,j=1,...,n$. Since this system is
maximal as a unimodular system, we obtain the negative answer on the
question 1). However, there are many other (maximal) unimodular
systems (see \cite{DG}) what gives many other classes of discrete
convexity.

The classes $\mathcal{P}t(\mathcal{R},\mathbb{Z})$ (as well as the class of integer
g-polymatroids) are stable under summation but not under intersection. Given an
unimodular system $\mathcal{R}$ one can construct another (dual) class of discrete
convexity which is stable under intersection (but not to summation). We show in
Theorem 3 that the theory becomes enough poor if to require stability  DC-class under
summation as well as under intersection.

It is worth to note that in sequel we develop the theory of discrete
convexity  not only for polytopes but for polyhedra as well. Because
of this we find more convenient to work with pure systems instead of
unimodular systems. Though, most interesting examples are related to
the latter ones.

The paper is organized as follows. In Section 2 we consider
several properties which one could want to require from a "good"
theory of discrete convexity. We find that all of them are in
essence equivalent. In Section 3 we introduce so called pure
systems and discuss their properties. In Section 4 we construct
classes of discrete convexity via the pure systems. Sections 5 and
6 are devoted to important particular case of pure systems, namely
to unimodular systems.
Each such a system enables us to construct a pair of (dual) DC-classes, one of which
is stable under summation and the other is stable under intersection, and these
classes contain ``many'' finite sets. In Section 7 we discuss an issue on defining
of polytopes from $\mathcal{P}t(\mathcal{R})$ by means of linear inequalities.

In a separate paper \cite{DK} we plan to develop corresponding
theory of discretely convex functions based on our theory of
discrete convexity. Let us note that particular cases of such a
theory relying on the DC class of g-polymatroids and its dual
DC-class was elaborated by Murota in series of papers, see, for
example \cite{M2,MSh}; in \cite{Kov2} was considered a class of
functions related to the DC-class dual to g-polymatroids (stable
under intersection), which was called later as L-convexity in
\cite{M2}.

Finally, we want to point out that recently the theory of discrete
convexity unexpectedly shown their importance in areas far from
discrete mathematics, such as in mathematical economics
\cite{DKM}, for a solution of the Horn problem \cite{steklov}, for
modules over discrete valuations rings \cite{izv}, in theory of
representation of groups.

\paragraph{Notations.}

In the sequel $M$ denotes a free Abelian group of finite
type\footnote{Of course, $M$ is isomorphic to $\mathbb Z^n$ for an
appropriate number $n$ but a general theory does not need to
distinguish a basis of the group.}. $V=M \otimes \mathbb
R\cong\mathbb R^n$ denotes the ambient vector space. Elements of $M$
are called {\em integer} points of $V$. Given a subset $P\subset V$,
we denote by $P(\mathbb Z)=P \cap M$ the set of integer points of
$P$.

$M^* =\mbox{Hom}(M,\mathbb Z)$ denotes the dual group, that is the
group of homomorphisms of Abelian groups $M \to \mathbb{Z}$. $V^*
=M^* \otimes \mathbb R$ is the dual vector space to $V$. For
$Q\subset V^*$, we put $Q(\mathbb Z)=Q \cap M^*$.

Let $X,Y$ be subsets of $V$. Then $X+Y=\{x+y,\, x \in X, \, y \in
Y\}$ denotes the (Minkowski) sum of $X$ and $Y$; $X-Y$ is understood
in a similar fashion. $\mbox{co}(X)$ denotes the convex hull of $X$
in $V$. $\mathbb Z(X)$ is the Abelian subgroup in $V$ generated by
$X$, that is the set of linear combinations of the form $\sum_x m_x
x$, where $x \in X$ and $m_x \in \mathbb Z$. $\mathbb R X$ denotes
the vector subspace generated by $X$.

\section{Discrete convexity: the basics}

The issue here is to characterize those subsets $X$ of a group $M$
($\cong{\mathbb Z}^n$) that we would be willing to call "convex"?
\medskip

{\bf Definition.} A subset $X \subset M$ is said to be {\em
pseudo-convex} if $X=\mbox{co}(X)({\mathbb Z})$ and $\mbox{co}(X)$
is a polyhedron.\medskip

Recall that a polyhedron is the intersection of some finite
collection of closed half-spaces of $V$. For example, a linear
sub-variety of $V$, or a polytope (the convex hull of some finite
subset in $V$) is a polyhedron. For more details about polyhedra,
see \cite{Grue} or \cite{R}.

We denote by $\mathcal PC$ the set of pseudo-convex sets.\medskip

{\bf Definition.} A polyhedron $P \subset V$ is {\em rational} if it
is given by a finite system of linear inequalities of the form
$p(v)\le a$ where $p\in M^*$ and $a\in \mathbb{Z}$. A polyhedron $P$ is called {\em
integer} if it is rational and if every (non-empty) face of $P$ contains an integer
point.\medskip

For example, a polytope is integer if and only if all its vertices
are integer points.\medskip

{\bf Proposition 1.} {\em  Suppose $X\subset M$. The following
assertions are equivalent:\smallskip

a)  $X$ is pseudo-convex;\smallskip

b)  $X=P({\mathbb Z})$ for some integer polyhedron $P \subset
V$;\smallskip

c) $X$ is the set of integer solutions of a finite system of linear
inequalities with integer coefficients.}\medskip

{\bf Proof.} The implication $a) \Rightarrow b)$ is almost obvious;
it suffices to take $P$ to be $\mbox{co}(X)$. The implication $b)
\Rightarrow c)$ is obvious. Finally, implication $c) \Rightarrow a)$
is precisely Meyer's theorem (see, for example, \cite{Sch}, Theorem
16.1). $\Box$\medskip

Denote by ${\mathcal IP}h$ the class of all integer polyhedra in
$V$. By Proposition 1, we have the natural bijection between the
classes ${\mathcal IP}h$ and ${\mathcal PC}$, which is given by
the mappings $P \mapsto P({\mathbb Z})$ and $X \mapsto
\mbox{co}(X)$. Both these classes are stable under integer
translations ($X\mapsto X+m$, $m\in\mathbb Z^n$), under the
reflection ($X\mapsto -X$), and under taking faces ($X\mapsto
X\cap F$, where $F$ is a face of the polyhedron $\mbox{co} (X)$).
Furthermore, the class ${\mathcal PC}$ is stable under
intersection and is not stable under summation, whereas the class
${\mathcal IP}h$ is stable under summation and is not stable under
intersection (the sum of two pseudo-convex sets needs not be
pseudo-convex, while the intersection of integer polyhedra need
not be integer).

Indeed, let us consider the following simple example in ${\mathbb
Z}^2$. Suppose $X=\{(0,0),(1,1)\}$ and $Y=\{(0,1),(1,0)\}$. Both $X$
and $Y$ are pseudo-convex. Despite that $X$ and $Y$ do not
intersect, they can not be separated by a linear functional (or a
hyperplane).

This example suggests that in order  to have the separation
property in theory of discrete convexity, we need to consider
narrower classes of subsets of $M$ than the class ${\mathcal PC}$.

We say that a class $\mathcal{K} \subset \mathcal{PC}$ is {\em
ample} if $\mathcal{K}$ is stable under  a) integer translations,
b) reflection, and c) faces. In the same way we understand
ampleness of a polyhedral class ${\mathcal P}\subset {\mathcal
IP}h$.\medskip

{\bf Proposition 2.} {\em  Let $\mathcal{K} \subset \mathcal{PC}$
be an ample class. The following four properties of $\mathcal{K}$
are equivalent:\smallskip

 $(Add)$ \ for every $X,Y\in \mathcal K$ the sets $X\pm Y$ are
pseudo-convex;\smallskip

$(Sep)$ \  if sets $X$ and $Y$ of $\mathcal K$ do not intersect,
then there exists (integer) linear functional $p:V\longrightarrow
\mathbb R$ such that $p(x)>p(y)$ for any $x\in X$, $y\in
Y$;\smallskip

$(Int)$ \ if sets $X$ and $Y$ of $\mathcal K$ do not intersect, then
the polyhedra   $co(X)$ and $co(Y)$ do not intersect as
well;\smallskip

$(Edm)$  \ for every $X,Y\in \mathcal K$ the polyhedron $co(X)\cap
co(Y)$ is integer. }\medskip

Proof. $(Add) \Rightarrow (Sep)$. \ If $X$ and $Y$ have an empty
intersection, then $0\notin X-Y$. Since the set $X-Y$ is
pseudo-convex, $0$ does  not belong to the polyhedron
$\mbox{co}(X-Y)=\mbox{co}(X)-\mbox{co}(Y)$. Hence there exists a
linear (integer) functional $p: V \to \mathbb R$ which is strictly
positive on $\mbox{co}(X-Y)$. Therefore $p(x)>p(y)$ for $x \in X$
and $y \in Y$.\medskip

$(Sep) \Rightarrow (Int)$. \ This one is obvious.\medskip

$(Int) \Rightarrow (Add)$. \ Let us show that $X-Y$ is
pseudo-convex. Since $\mbox{co}(X-Y)= \mbox{co}(X)-\mbox{co}(Y)$ is
a polyhedron, we need to prove that $X-Y=\mbox{co}(X-Y)\cap M$.
Suppose the integer point $m$ lies in
$\mbox{co}(X-Y)=\mbox{co}(X)-\mbox{co}(Y)$. Then the polyhedra
$\mbox{co}(X)$ and $m + \mbox{co}(Y)=\mbox{co}(m+Y)$ intersect.
Applying $(Int)$ to the sets $X$ and $m+Y$, we see that these sets
also intersect, that is $m \in X-Y$.\medskip

$(Edm) \Rightarrow (Int)$. \ This implication is obvious.\medskip

$(Int) \Rightarrow (Edm)$. \ Suppose $X,Y \in \mathcal K$,
$P=\mbox{co}(X)$, $Q=\mbox{co}(Y)$. We need to show that $P \cap Q$
is an integer polyhedron. Obviously $P \cap Q$ is rational.
Therefore we need to establish that every (non-empty) face of $P\cap
Q$ contains an integer point. We assume here, without loss of
generality, that the face is minimal.

Suppose $F$ is a minimal (non-empty) face of the polyhedron $P\cap Q$.
Let $P'$ (resp. $Q'$) be a minimal face of $P$ (resp. $Q$) which
contains $F$. We claim that  $F= P'\cap Q'$.

Projecting $V$ along $F$, we may suppose additionally that
$F$ is of dimension $0$. That is $F$ consists of a single point, which
is a vertex of $P \cap Q$. Suppose, on the contrary, that $P'\cap Q'$
contains some other point $a$. Since the point $F$ is relatively
interior both in $P'$ and in $Q'$, then $F$ is an interior point of
some segment $[a,b]$, lying in both $P'$ and $Q'$. But in such a case
the segment $[a,b] \subset P'\cap Q'\subset P\cap Q$, and $F$ can not
be a vertex of $P\cap Q$. Contradiction.

Thus, $F= P'\cap Q'$. Since our class $\mathcal K$ is stable under
faces, the sets $P'({\mathbb Z})$ and $Q'({\mathbb Z})$ belong to
$\mathcal{K}$. The property {\em (Int)} implies that the sets
$P'({\mathbb Z})$ and $Q'({\mathbb Z})$ intersect. Because of
this, $F$ is an integer singleton. $\Box$\medskip

{\bf Definition.} An ample class $\mathcal K\subset \mathcal PC$ is
a {\em class of discrete convexity} (or a {\em DC-class}) if it
possesses anyone of the properties from Proposition 2.\medskip

On the language of integer polyhedra, the definition of discrete
convexity is formulated as follows. A class $\mathcal{P}$ of integer
polyhedra is a {\em polyhedral class of discrete convexity} if it is
ample and the following variant of the Edmonds' condition
holds:\medskip

{$ (Edm')$} The intersection of any two polyhedra from ${\mathcal
P}$ is an integer polyhedron (not necessarily in ${\mathcal
P}$).\medskip

According to Proposition 2, the equivalent requirement is:\medskip

{$ (Add')$} \ \ \ $(P+Q)({\mathbb Z})=P({\mathbb Z})+Q({\mathbb Z})$
for every $P,Q \in \mathcal{P}$.\medskip

Let us give a few examples of DC-classes.\medskip

{\bf Example 1. One-dimensional case.} Let $M\cong {\mathbb Z}$.
Then the class $\mathcal{PC}$ of all pseudo-convex sets is a {\em
DC}-class. This is not the case in higher dimensions of course.
$\Box$\medskip

The class of integer rectangles in the plane $\mathbb R^2$ is a {\em
DC}-class. More generally, if ${\mathcal{K}}_1$ and ${\mathcal K}_2$
are {\em DC}-classes in the groups $M_1$ and $M_2$, respectively,
then the class of sets of the form $X_1\times X_2$ with $X_i \in
{\mathcal K}_i$, $i=1, 2$, is a {\em DC}-class in $M_1\times M_2$ as
well.\medskip

{\bf Example 2. Hexagons}. Let us consider a more interesting
class ${\mathcal H}$ of polyhedra in $\mathbb R^2$. It consists of
polyhedra defined by the inequalities $a_1 \le x_1 \le b_1$,
$a_2\le x_2 \le b_2$, $c \le x_1+x_2\le d$, where $a_1$, $a_2$,
$b_1$, $b_2$, $c$ and $d$ are integers. It is easy to check that
this hexagon (generally speaking, this hexagon can be degenerated
to a polyhedron with smaller number of edges) has integer
vertices. Obviously, ${\mathcal H}$ is stable under integer
translations, reflection and faces. Since the intersection of
hexagons yields a hexagon, we conclude that ${\mathcal H}$ is a
polyhedral {\em DC}-class. $\Box$\medskip

{\bf Example 3. Base polyhedra.} This is one of the possible
high-dimensional generalizations of Example 2.  Let $N$ be a finite
set, and $V=(\mathbb R^N)^*$. We interpret elements of $V$ as
measures on the set $N$.  Recall, that a function $b: 2^N \to
\mathbb R \cup\{+\infty\}$ is called {\em submodular} if for any
$S$, $T \subset N$, the following inequality holds
$$
b(S)+ b(T) \ge b(S \cup T)+ b(S\cap T).
$$
The elements of $V$ can be viewed as modular functions, i.e.,
functions which fulfill the above-written definition of
submodularity with equality.

A {\em  base polyhedron} is a polyhedron of the following form
$$
B(b)=\{x \in V \ | \ x(S) \le b(S), \, S \subset N, \,
\mbox{ and  } x(N)=b(N)\},
$$
where $b$ is a submodular function. Obviously, the class
${\mathcal B}$, which consists of base polyhedra with
integer-valued $b$, is stable under integer translations and under
reflection. One can show that it is stable under faces, and hence,
each base polyhedron has integer vertices. The well-known theorem
by Edmonds \cite{Ed} ensures that the property ({\em Edm})
obtains, and thus ${\mathcal B}$ is a polyhedral {\em DC}-class.
The reader will find details of the proofs of these properties of
base polyhedra in \cite{FT}, or see our Example 13. $\Box$\medskip

{\bf Example 4.} Here we give another high-dimensional
generalization of Example 2. Let $N$ be a finite set, and let
$V=\mathbb R^N$ be the space of real-valued functions on $N$.
Consider the class ${\mathcal L}$ of polyhedra in $V$, given by the
inequalities of the form $a_i \le x(i) \le b_i$ and $a_{ij} \le
x(i)-x(j) \le b_{ij}$, where $i,j \in N$, and all $a$'s and $b$'s
are integers. We claim that these polyhedra are integer. Indeed,
their vertices are given by equalities of the form $x(i) =c_i$ and
$x(i) -x(j)=c_{ij}$ where $c$'s are integers. It is clear that $x$
is an integer point.

Thus, the class ${\mathcal L}$ consists of integer polytopes.
Since it is stable under intersection, the axiom $(Edm')$ is
satisfied automatically, and ${\mathcal L}$ is a polyhedral {\em
DC}-class. $\Box$\medskip

We give a general construction of {\em DC}-classes in Section
4.\smallskip

In the classical context, convexity is preserved under summation
and under intersection. It would be natural therefore to require
these properties for the discrete set-up. For example, both the
classes of segments and hexagons and their products possess these
properties. Moreover (see  Theorem 2), these cases exhaust {\em
DC}-classes, stable under both summation and intersection. The
class $\mathcal B$ described in Example 3 is stable under
summation, but not under intersection (if $|N|>3$). Similarly, the
class ${\mathcal L}$ described in Example 4 is stable under
intersection, but not under summation (if $|N|>2$). Therefore,
when we consider classes stable under summation and classes stable
under intersection separately, more interesting theory of discrete
convexity is obtained.\medskip

{\bf Definition.} An ample class ${\mathcal K}$ of pseudo-convex
sets is called an {\em S-class} if $X+Y \in {\mathcal K}$ for any
$X,Y \in {\mathcal K}$.\medskip

In particular, $X-Y \in \mathcal{PC}$ for any $X,Y \in
\mathcal{K}$, and, thus, any {\em S}-class is a {\em DC}-class.
However in order to characterize polyhedral {\em S}-classes, we
have to require both that the class be stable under summation and
the axiom $(Add')$ be satisfied. Note that the intersection of two
polyhedra of a polyhedral {\em S}-class is an integer polyhedron,
but need not be a polyhedron of this class.\medskip

{\bf Definition.} An ample class $\mathcal{P}$ of integer polyhedra
is called a {\em polyhedral} I-{\em class} if $P \cap Q\in
\mathcal{P}$ for any $P,Q\in \mathcal{P}$.\medskip

Again any {\em I}-class is a {\em DC}-class, since the axiom
$(Edm')$ holds. Let $P$ and $Q$ be polyhedra in an {\em I}-class,
then $P({\mathbb Z})+Q({\mathbb Z})$ is a pseudo-convex set, though
$P+Q$ need not be a polyhedron of this class.

\section{Pure systems}

Linear subspaces are the simplest polyhedra. For a (rational) vector
subspace $F \subset V$ the set $S=F(\mathbb{Z})$ of all integer
points of $F$ is an Abelian subgroup of $M$. Such subgroups of $M$
are called {\em pure}. Let us collect some properties of pure
subgroups (of $M$) in the following simple\medskip

{\bf Lemma 1.} {\em Let $S$ be a subgroup of a free Abelian group of
finite type $M$. The following assertions are equivalent:

1) $S$ is a pure subgroup;

2) $S$ is a pseudo-convex subset of $M$;

3) the factor-group $M/S$ is torsion-free;

4) the factor-group $M/S$ is a free Abelian group.} $\Box$\medskip

In fact, the factor-group $M'/f^{-1}(S)$ is imbedded in the
torsion-free group $M/S$ and, therefore, has no torsion.
$\Box$\medskip

In general, the sum of pure subgroups of $M$ need not be a pure subgroup of $M$.
For  example, if $M=\mathbb{Z}^2$, $S=\mathbb{Z}(1,1)$, $S'=\mathbb{Z}(1,-1)$ then
the group $S+S'$ has the index 2 in $M$.\medskip

{\bf Definition.} Pure subgroups $S$ and $S'$ of $M$ are called {\em
mutually pure} if the sum $S+S'$ is a pure subgroup of $M$. Two
(rational) linear subspaces $L$ and $L'$ of $V$ are {\em mutually
pure} if the subgroups $L(\mathbb{Z})$ and $L'(\mathbb{Z})$ are mutually
pure.\medskip

There is the following criterion of the mutual purity.\medskip

{\bf Lemma 2.} {\em Let $S_1$ and $S_2$ be two pure subgroups of
$M$. They are mutually pure if and only if the image of natural
homomorphism $S_1 \to M/S_2$ is pure.}\medskip

In fact, the factor group $(M/S_2)/(\text{Im}(S_1))$ is canonically
isomorphic to $M/(S_1+S_2)$. $\Box$\medskip

Pure subgroups naturally come in play in the study of DC-classes. Suppose we have a
pseudo-convex subset $X$ in $M$. Then we can consider the linear subspace
$Tan(X):=\mathbb{R}(X-X)$ in $V$ (the "tangent space" of $X$) and the subgroup
$S=\mathbb{Z}(X-X)$ in $M$. Of course, $S \subset Tan(X)(\mathbb{Z})$, and  in the
general case this inclusion is proper. Hence, in the general case,
$S$ needs not  be a pure subgroup of $M$. Nevertheless, there is
an instance when we can guarantee the purity of $S$.

For a natural number $n$ and $X\subset M$, we denote by $[n]X$ the
sum of $n$ copies of $X$; for example, $[2]X=X+X$.\medskip

{\bf Proposition 3.} {\em Let $X \subset M$. Suppose that $[n]X$
is a pseudo-convex set for every $n=1,\ldots $. Then the subgroup
$\mathbb{Z}(X-X)$ is pure.}\medskip

Proof. Changing $X$ by $[n]X$ for an appropriate large $n$, one can assume that $X$
contain a point $a$ which belongs to the relative interiority of co$(X)$. Changing
$X$ by $X-a$, one can assume that $0$ belongs to the relative interiority of co$(X)$.
In that case  $\mathbb{Z}(X-X)=\cup_{n\ge 1}[n]X$. It remains to note that an
increasing union of pseudo-convex sets is a pseudo-convex set. $\Box$\medskip

Given an ample class $\mathcal{K}$ of pseudo-convex sets, we can
associate to it the following system $\mathcal{U}(\mathcal{K})$ of
linear subspaces in $V$ (the {\em homogenization} of
$\mathcal{K}$). Namely,
$$
{\mathcal U}({\mathcal K})=\{Tan(X), \ X\in {\mathcal K}\}.
$$
Similarly we define the system of vector subspaces
$\mathcal{U}(\mathcal{P})$ for an ample polyhedral class
$\mathcal{P}$.\medskip

{\bf Definition.} A collection $\mathcal{U}$ of linear subspaces
in $V$ is called  a {\em pure system} if every $F, G \in
\mathcal{U}$ are mutually pure subspaces. Elements of a pure
system are called {\em flats}.\medskip

The homogenization of DC-classes produces pure systems. Say
that an ample class $\mathcal{P}$ of integer polyhedra is {\em very ample} if it
contains the polyhedron $nP$ with any integer $n$ and any polyhedron $P \in
\mathcal{P}$.\medskip

{\bf Proposition 4.} {\em Let $\mathcal{P}$ be a very ample DC-class
$\mathcal{P}$ of integer polyhedra. Then $\mathcal{U}(\mathcal{P})$
is a pure system.}\medskip

Proof. Let $F=Tan(P)$ and $G=Tan(Q)$, where $P,Q \in \mathcal{P}$.
We have to show that the subgroup $F(\mathbb{Z})+G(\mathbb{Z})$ is
pure. Of course, this subgroup contains the subgroup
$\mathbb{Z}((A+B)-(A+B))$, where $A=P(\mathbb{Z})$ and
$B=Q(\mathbb{Z})$. According to Proposition 3, it suffices to
check that the set $[n](A+B)=[n]A+[n]B$ is pseudo-convex for any
$n=1,\ldots $.

Since the class $\mathcal{P}$ is discretely convex, the set $A+A$
is pseudo-convex and coincides with $2P(\mathbb{Z})$. Similarly,
for any $n$, $[n]A=(nP)(\mathbb{Z})$ as well as
$[n]B=nQ(\mathbb{Z})$. At last,
$[n]A+[n]B=nP(\mathbb{Z})+nQ(\mathbb{Z})=(nP+nQ)(\mathbb{ Z})$ is
a pseudo-convex set since $nP$ and $nQ$ belong to $\mathcal{P}$.
$\Box$\medskip

In the next Section we show how to dehomogenize pure
systems.\smallskip

A pure system $\mathcal{U}$ is said to be a {\em pure S-system}
(correspondingly, a {\em pure I-system}) if $F+G$
(correspondingly, $F\cap G$) belongs to $\mathcal{U}$ for any $F,G
\in \mathcal{U}$. It is clear that the homogenization of an
S-class is a pure S-system, and the homogenization of an I-class
is a pure I-system.\medskip

Let us illustrate the homogenization procedure on the class
${\mathcal B}$ of base polyhedra from Example 3.\medskip

{\bf Example 5. The homogenization of base polyhedra.} Recall, that
here $V=(\mathbb R^N)^*$ is the space of measures on a finite set
$N$. Let $B(b)$ be the base polyhedron defined by a submodular
function $b:2^N \rightarrow \mathbb R \cup \{+\infty\}$. We are
going to show how the corresponding tangent space $Tan(B(b))$ looks
like. Here we can assume that $B(b)$ is a symmetric (with respect to
the origin $0$) base polyhedron. This means that $b(S)=b(N\setminus
S)$; in particular, $b(N)=0$. It is clear, that $nB(b)=B(nb)$. Therefore the tangent
space $Tan(B(b))$ is the base polyhedron $B(\infty b)$, that is given by the
following list of equations $$ x(S)=0,\,\ S\in{\mathcal F}(b), $$ where
$\mathcal{F}(b)=\{S\subset N, \ b(S)=0$\}. Obviously, $\emptyset, N \in
\mathcal{F}(b)$. The symmetry of $B(b)$ implies that $N\setminus S \in
\mathcal{F}(b)$ with any $S \in \mathcal{F}(b)$. Submodularity of $b$ implies that
$S\cup T$ and $S\cap T$ belong to $\mathcal{F}(b)$ with any $S,T \in
\mathcal{F}(b)$. Thus, $\mathcal{F}(b)$ is a Boolean subalgebra of
$2^N$.

We see that to give a flat of $\mathcal{U}(\mathcal{B})$ is the
same as to give a Boolean subalgebra of $2^n$, or is the same as
to give an equivalence relation $\approx$ on $N$. The
corresponding flat $F(\approx)$ consists of measures $x\in V$ such
that $x(S)=0$ for each equivalence class $S$ of the relation
$\approx$. The codimension of this flat $F(\approx)$ is equal to
the number of equivalence classes of  $\approx$.

Let us consider, for instance, one-dimensional flats. These flats
correspond to those equivalence relations which possess a single
class of equivalence of cardinality $2$, whereas all others
classes are of cardinality $1$. For example the one-dimensional
flat $\mathbb R(e_i-e_j)$  corresponds to the equivalence relation
whose $2$-element class of equivalence is $\{i,j\}$. Here $(e_i)$,
$i\in N$, denote the Dirac measure at the point $i\in N$.

Similarly, flats of codimension $1$ correspond to dichotomous
equivalence relations (i.e., relations with only two equivalence
classes, say $T$ and $N\setminus T$).

We denote $\mathcal{U}(\mathbb A(N))$ this pure system.
$\Box$\medskip

Let us return to general pure systems. There holds the following
finiteness property.
\medskip

{\bf Proposition 5.} {\em Any pure system is a finite set.}\medskip

Proof. Let $\mathcal{U}$ be a pure system  of pure subgroups in
$M$. Let $\mathbb{F}_2$ be the 2-elements field. For any pure subgroup $S$ in $M$ we
can consider the corresponding $\mathbb{F}_2$-vector subspace $S\otimes
\mathbb{F}_2$ in the vector space $M\otimes \mathbb{F}_2$. It is clear that the
dimension of $S\otimes \mathbb{F}_2$ is equal to the rank of $S$ (that is the
dimension of $S\otimes \mathbb{R}$).

We assert that for different $S,S' \in \mathcal{U}$ their images
$S\otimes \mathbb{F}_2$ and $S'\otimes \mathbb{F}_2$ are also different. Suppose
that $S\otimes \mathbb{F}_2=S'\otimes \mathbf{F}_2$. Then $(S+S')\otimes
\mathbb{F}_2=(S\otimes \mathbb{F}_2)+(S'\otimes \mathbb{F}_2)=S\otimes
\mathbb{F}_2$. Since $S+S'$ is pure then the rank of $S+S'$ is equal to the rank of
$S$ (and is equal to the rank of $S'$). Therefore $S=S+S'=S'$. $\Box$\medskip

{\bf Dualization.}  Now we discuss a construction of dual (or
orthogonal) pure system. For a vector subspace $L$ in $V$, let
$L^{\bot}$ denote the orthogonal vector subspace in the dual vector
space $V^*$, that is
$$
L^{\bot}=\{ p\in V^*, p(v)=0 \text{ for any } v\in L\}.
$$

{\bf Theorem 1}. {\em If $L$ and $L'$ are mutually pure subspaces in $V$ then
$L^{\bot}$ and $L'^{\bot}$ are mutually pure subspaces in $V^*$.}\medskip

For proving this theorem, it is convenient to use a notion of a pure
homomorphism. Let $M$ and $N$ be free Abelian groups of finite type. Let us say that
a homomorphism $f:M\rightarrow N$ is {\em pure} if the factor-group $N/f(M)$ is a
free (or torsion-free) Abelian group. This means, of course, that $f(M)$ is a pure
subgroup in $N$.\medskip

{\bf Lemma 3.} {\em A homomorphism $f:M\rightarrow N$ is pure if and
only if the dual homomorphism $f^*:N^*\rightarrow M^*$ is
pure.}\medskip

Proof. Let us consider the canonical decomposition of the
homomorphism  $f:M\rightarrow N$ in two exact sequences
$$
0\rightarrow K\rightarrow M\rightarrow H\rightarrow 0, \quad\mbox{
and }\quad 0\rightarrow H\rightarrow N\rightarrow C\rightarrow 0.
$$
Since $f$ is pure, $C$ is a free Abelian group. The group $H$ is
free as a subgroup of the free group $N$. Therefore, both
sequences are split. Hence the dual sequences
$$
0\rightarrow C^*\rightarrow N^*\rightarrow H^* \rightarrow 0, \quad
0\rightarrow H^*\rightarrow M^*\rightarrow K^* \rightarrow 0
$$
are exact (where $X^*=\mbox{Hom} (X,\mathbb Z)$). Since $K^*$ is
free, we obtain that $f^*:N^*\rightarrow M^*$ is pure.
$\Box$\medskip

Proof of Theorem 1. Let $S=L(\mathbb{Z})$, and similarly $S'=L'(\mathbb{Z})$. It is
obvious that $L^{\bot}(\mathbb{Z})$ is equal to $S^{\bot}=\{ p\in M^*, p(s)=0 \
\forall \ s\in S\}$. That is that $S^{\bot}$ is the kernel of the canonical
projection $M^* \to S^*$ being dual to the inclusion $S \to M$. It is clear from
this, that $(S^{\bot})^*$ can be identified with $M/S$.

We have to show that the subgroup $S^{\bot} + S'^{\bot}$ is pure
in $M^*$. That is, by Lemma 2,  that the canonical homomorphism
$S^{\bot} \to M^*/S'^{\bot}$ is pure. By Lemma 3, it suffices to
check that the dual homomorphism $(M^*/S'^{\bot})^* \to
(S^{\bot})^*$ is pure. The latter homomorphism can be identified
with the canonical homomorphism $S'^{\bot} \to M/S$. But this
homomorphism is pure because $S$ and $S'$ are mutually pure
subgroups. $\Box$\medskip

{\bf Corollary.} {\em Let $\mathcal{U}$ be a pure system in $V$.
Then the collection $\mathcal{U}^{\bot}:=\{L^{\bot}, L\in
\mathcal{U}\}$ is a pure system in $V^*$.} $\Box$

\section{Construction of DC-classes}
In the previous section, we constructed pure systems via the
homogenization of (very ample) DC-classes. Here we shall go in the
opposite direction.

Let $\mathcal{U}$ be a pure system in $V$. If we consider all
integer translations of flats of $\mathcal{U}$, we obtain a
polyhedral DC-class. However, this class is of a little interest.
For instance, it contains no polytopes (except, may be,
0-dimensional ones). Below we define a more interesting (maximal)
DC-class ${\mathcal P}h({\mathcal U},\mathbb Z)$ of integer
polyhedra associated to a given pure system ${\mathcal U}$.\medskip

{\bf Definition.} Let ${\mathcal U}$ be a collection of (rational)
vector subspace in $V$. A polyhedron $P$ is said to be ${\mathcal
U}$-{\em convex} (or ${\mathcal U}$-{\em polyhedron}) if, for any
face $F$ of $P$, the tangent space $Tan(F)=\mathbb R (F-F))$
belongs to ${\mathcal U}$.\medskip

Let ${\mathcal P}h({\mathcal U})$ be the set of ${\mathcal
U}$-polyhedra, and let  ${\mathcal P}h({\mathcal U},\mathbb Z)$ be
the set of integer ${\mathcal U}$-polyhedra. Note that the class
${\mathcal P}h({\mathcal U},\mathbb Z)$ is stable under integer
translations, reflection and faces. In other words, it is an ample
(and even very ample) class of integer polyhedra. The
homogenization of ${\mathcal P}h({\mathcal U},\mathbb Z)$ brings
us back to ${\mathcal U}$.

The following result will be used in the sequel.\medskip

{\bf Proposition 6.} {\em Let $P\in\mathcal{P}h(\mathcal{U}, \mathbb{Z})$, and $L$
be an integer vector subspace in $V$. Suppose that $L$ is mutually pure with any
subspace of $\mathcal{U}$. Them the intersection $P\cap L$ is an integer
polyhedron.}\medskip

Proof. Let $\gamma $ be a minimal face of $P\cap L$; we have to
show that $\gamma$ is an integer polyhedron. In fact, $\gamma$ is
an affine subspaces in $V$ because it has no faces. Changing $P$
by its minimal face containing $\gamma$ we may assume that
$\gamma=P\cap L$. Now, if we replace $P$ by its affine span
aff$(P)$, then we would have aff$(P)\cap L=\gamma$. But aff$(P)$
is an integer translation of the linear subspace
$\mathbb{R}(P-P)$. Therefore we can assume that $P$ is an integer translation of a
linear subspace $L'$ in $V$, $P=L'+m$ for some $m\in M$.

Now we can repeat the reasoning from Proposition 2. If $L$ and
$L'+m$ do not intersect, then the assertion is obviously true. Let
$x\in L\cap (L'+m)$, that is $x\in L$ and $x=x'+m$, $x'\in L'$.
Then $m=x-x'$ is an integer point of $L-L'$. Since $L$ and $L'$
are mutually pure, $m\in L(\mathbb{Z})-L'(\mathbb{Z})$. That is, there exists an
integer point $l\in L$ such that $l+m\in L'$. $\Box$\medskip

Now we show that if $\mathcal{U}$ is a pure system, then
$\mathcal{P}h(\mathcal{U},\mathbb{Z})$ is a DC-class.\medskip

{\bf Theorem 2.} {\em A class $\mathcal{P}h(\mathcal{U},\mathbb{Z})$ is a DC-class of
integer polyhedra if and only if $\mathcal{U}$ is a pure system.}\medskip

Proof. Since $\mathcal{P}h(\mathcal{U},\mathbb{Z})$ is a very ample,
the "only if" part of Theorem war proven in Proposition 4. Let us
prove the "if" part.

More precisely, we shall show that the intersection of two
polyhedra from the class $\mathcal{P}h(\mathcal{U},\mathbb{Z})$ is
an integer polyhedron. For this, we use a trick known in Algebraic
Geometry as "reduction to the diagonal". Namely, we replace the
intersection of two polyhedra $P$ and $Q$ by the intersection of
their direct product $P\times Q$ with the linear subspace $\Delta$
being the diagonal in $V\times V$.

Let us consider in $V\times V$ the following system
$\mathcal{U}\times \mathcal{U}$ of subspaces of the form $L\times
L'$ where $L,L'\in \mathcal{U}$. Obviously, $P\times Q$ is
$\mathcal{U}\times \mathcal{U}$-polyhedron. The intersection
$P\times Q$ with the diagonal $\Delta$ consists of points of the
form $(v,v)$ such that $v$ belongs to $P$ and to $Q$. Therefore to
prove that $P\cap Q$ is an integer polyhedron in $V$ is the same
as to prove that $(P\times Q)\cap \Delta$ is an integer polyhedron
in $V\times V$.

By virtue of Proposition 7, it suffices to show that the diagonal
$\Delta$ is mutually pure with any subspace $L\times L'$ where
$L$, $L'\in \mathcal{U}$. But this is equivalent to the mutual
purity of the subspaces $L$ and $L'$. The latter property holds by
the definition of pure systems. $\Box$\medskip

{\bf Remark.} Using the above arguments we obtain the following
more general result. Let $\mathcal{U}$ and $\mathcal{U'}$ be two
systems of subspaces in $V$. Suppose that for every $L\in
\mathcal{U}$ and $L'\in \mathcal{U'}$ the subspaces $L$ and $L'$
are mutually pure. Then the intersection of any integer
$\mathcal{U}$-polyhedron with any integer
$\mathcal{U'}$-polyhedron is an integer polyhedron.

Of course, if a pure system $\mathcal{U}$ is stable under
summation (intersection) then the corresponding class
$\mathcal{P}h(\mathcal{U},\mathbb{Z})$ is an S-class (I-class).

\section{Unimodular systems}
We have shown above, that pure systems play a crucial role in the
description and construction of $DC$-classes (of integer polyhedra
in $V$ or of pseudo-convex subsets in $M$). The corresponding
$DC$-classes contain, for example, all integer translations of
flats. However, if we want that a $DC$-class contains {\em
polytopes}, we have to provide that the corresponding pure system
has "sufficiently many" one-dimensional flats. This means that every
flat of our system is generated (as a vector subspace) by
one-dimensional flats. Here we explain how to construct pure
S-systems (and the dual pure I-systems, see the next Section) by
means of unimodular systems.\medskip

{\bf Definition.} A subset ${\mathcal R} \subset M$ is called {\em
unimodular} if, for any subset $B \subset {\mathcal R}$ the subgroup
${\mathbb Z} B \subset M$ is pure. A {\em unimodular system} is a
pair $(M, {\mathcal R})$ where ${\mathcal  R}$ is a unimodular set
in $M$. Non-zero elements of $\mathcal{R}$ are called {\em
roots}.\medskip

We call {\em flats} (or ${\mathcal R}$-{\em flats}) subspaces
$\mathbb{R} B$, where $B \subset {\mathcal R}$. It is obvious that
the set $\mathcal{U}(\mathcal{R})$ of all $\mathcal{R}$-flats is a
pure S-system.

Unimodular systems are closely related to totally unimodular
matrices, that is matrices whose minors are equal to  $0$ or $\pm
1$. Suppose that a unimodular set ${\mathcal R}$ is of full
dimension, or, equivalently, spans $V$. If we pick a basis $B
\subset {\mathcal R}$ and represent vectors of ${\mathcal R}$ as
linear combinations of the basis vectors, then the matrix of
coefficients is totally unimodular. In particular, the
coefficients of this matrix are either $0$ or $\pm 1$, which
proves finiteness of any unimodular set. Conversely, columns of a
totally unimodular $n\times m$ matrix yield a unimodular set in
$\mathbb Z^n$. Thus unimodular systems are nothing but
coordinate-free representations of totally unimodular matrices.
The reader might find many other characterizations of totally
unimodular matrices in \cite{Sch}.

Consider some important examples of unimodular systems.\medskip

{\bf Example 6}. In Example 5, we introduced the pure system
$\mathbb A(N)$, which is spanned by one-dimensional flats $\mathbb
Z(e_i-e_j)$, $i,j\in N$. Therefore, the set of vectors $e_i-e_j$,
$i,j \in N$, is a unimodular set in $(\mathbb Z^{N})^*$. Let us
denote this system as well by ${\mathbb A}(N)$. Note that it is not
of full dimension, since it spans the subspace $\{x, x(N)=0\}$,
which is orthogonal to the vector ${\bf 1}_N \in \mathbb R^N$. We
shall show in Section 7 that the class  ${\mathcal P}h({\mathbb
A(N)})$ coincides with the class of base polyhedra ${\mathcal B}$
from Example 3.

If we project the set $\mathbb A(N \cup \{0\}))$ along the axis
$\mathbb R e_0$ onto the space $(\mathbb R^N)^*$, we obtain the
full-dimensional unimodular system consisting of the vectors $\pm
e_i$ and $e_i-e_j$, $i,j \in N$, in $(\mathbb Z^N)^*$. Of course, we
could construct this system simply by adding the basic system $(\pm
e_i, i \in N)$ to the system  $\mathbb A(N)$. We denote this system
by $\mathbb A_N$. We shall show that $\mathbb A_N$-polyhedra are
precisely generalized polymatroids.

Sub-systems ${\mathcal R} \subset \mathbb A_N$ (more precisely,
symmetrical sub-systems, which contain $0$ and $-r$ for any $r \in
{\mathcal R}$) are called {\em graphic} unimodular systems.\medskip

{\bf Example 7}. To any graph $G$ one can associate another
unimodular system, the so called  {\em cographic} unimodular system
$\mathbb D(G)$. It is located in the cohomology group $H^1(G,\mathbb
Z)$ of the graph $G$ and consists of the cohomology classes $\pm
[e]$, corresponding to oriented edges of the graph $G$. The proof of
the unimodularity of the system $\mathbb D(G)$ is based on the fact
that this system is, in some (matroidal) sense, dual to the graphic
system associated with $G$.

Cubic (or 3-valent) graphs gives the most interesting examples of
cographic systems. The simplest example of such a graph is the
complete graph $K_4$ with $4$ vertices. The corresponding system
$\mathbb D (K_4)$ is isomorphic to ${\mathbb A}_3$. The bipartite
graph $K_{3,3}$ yields a more interesting example. The system
$\mathbb D (K_{3,3})$ consists of the following $19$ vectors in
$\mathbb R^4$: $\{0$, $ \pm e_i$, $i=1,\ldots, 4$, $\pm (e_1+e_2)$,
$\pm (e_2+e_3)$, $\pm (e_3+e_4)$, $\pm (e_4+e_1)$, $\pm
(e_1+e_2+e_3+e_4)\}$.

One can check that $\mathbb D (K_{3,3})$ is not a graphic system.
$\Box$\medskip

{\bf Example 8}. There is an exceptional unimodular system $\mathbb
E_5$ in dimension $5$ which is neither graphic no cographic. It
consists of the following $21$ vectors: $0$, $\pm e_i$, $i=1,\ldots
,5$, $\pm (e_1-e_2+e_3)$, $\pm (e_2-e_3+e_4)$, $\pm (e_3-e_4+e_5)$,
$\pm (e_4-e_5+e_1)$, $\pm (e_5-e_1+e_2)\}$. $\Box$\medskip

According to the Seymour theorem \cite{Sey}, every unimodular system
can be constructed via graphic systems, cographic systems, and the
system $\mathbb E_5$.\smallskip

Let $(M,\mathcal{R})$ and $(M',\mathcal{R'})$ be unimodular systems.
A homomorphism of Abelian groups $\varphi: M \to M'$ is called a
{\em morphism of unimodular systems} if $\varphi(\mathcal{R})
\subset \mathcal{R'}$. For example, if $\varphi$ is the projection
of $M$ onto $M'=M/\mathbb{Z}r$, where $r\in \mathcal{R}$, then
$\varphi (\mathcal{R})$ is a unimodular set in $M'$. The {\em direct
sum} of unimodular systems $(M,\mathcal{R})$ and $(M',\mathcal{R'})$
is a unimodular system $(M\oplus M', \mathcal{R}\oplus
\mathcal{R'})$, where $\mathcal{R}\oplus
\mathcal{R'}=\mathcal{R}\cup \mathcal{R'}$.

The following theorem characterizes unimodular systems ${\mathcal
R}$ whose pure systems ${\mathcal U}({\mathcal R})$ are stable
under intersection. For such a system the corresponding DC-class
$\mathcal{P}(\mathcal{U})$ is simultaneously  S-class and
I-class.\medskip

{\bf Theorem 3.} {\em Let ${\mathcal R}$ be a unimodular set such that the pure
system ${\mathcal U}({\mathcal R})$ is stable under intersection. Then ${\mathcal R}$
is the direct sum of copies of $\mathbb A_1$ and $\mathbb A_2$.}\medskip

Proof. The proof is by induction on the dimension of unimodular
systems. The assertion is obvious in dimensions $1$ and $2$.

Consider first of all the case of dimension $3$. Assume ${\mathcal
R}$ contains a flat $S$ isomorphic to $\mathbb A_2$. Denote by
$e_1$, $e_2$ and $e_1+e_2$ the vectors of ${\mathcal R}\cap S$. We
claim that there is at most one more vector of ${\mathcal R}$ (up
to collinearity). Suppose there are two non-collinear vectors.
Clearly we may denote them by $e_3$ and $e_1+e_3$. Then, since
${\mathcal U}({\mathcal R})$ is stable under intersection,
$e_2-e_3$ and $e_1+e_2+e_3$ belong to ${\mathcal R}$. But this
contradicts unimodularity of ${\mathcal R}$, and the claim is
proven. Therefore, ${\mathcal R}$ is isomorphic to $\mathbb
A_1\oplus\mathbb A_2$.

One can similarly check that if ${\mathcal R}$ does not contain
flats isomorphic to $\mathbb A_2$, then ${\mathcal R}$ is isomorphic
to $\mathbb A_1\oplus\mathbb A_1\oplus\mathbb A_1$. Thus, in the
$3$-dimensional case, the proposition  is verified.

General case. Let ${\mathcal U}({\mathcal R})$ contain a flat $S$
isomorphic $\mathbb A_2$. This means that $S$ is a plane of $V$ such
that ${\mathcal R}\cap S\cong \mathbb A_2$. We will show that there
exists a flat $T$ of codimension $2$ in $V$ such that
\begin{equation}\label{deco-1}
{\mathcal R}=({\mathcal R}\cap S)\cup ({\mathcal R}\cap T).
\end{equation}
By induction ${\mathcal R}\cap T$ is equal to the sum of copies
$\mathbb A_1$ and $\mathbb A_2$, and we have ${\mathcal R}\cap
S\cong \mathbb A_2$, so if (\ref{deco-1}) is true, the proposition
is also true.

Pick a flat $T$ of ${\mathcal U}({\mathcal R})$ of codimension $2$
(in $V$) such that $T\cap S=0$. Obviously such a flat
exists.\medskip

{\em Claim}. ${\mathcal R}\subset S\cup T$.\medskip

Let us consider the projection $\pi:V\rightarrow S$ which has $T$
as the kernel (the projection along $T$). Then $\pi ({\mathcal
R})$ is a unimodular system of $S$ which contains ${\mathcal
R}\cap S$. Because ${\mathcal R}\cap S\cong \mathbb A_2$ and the
$\mathbb A_2$ is a maximal unimodular system, any vector
$r\in{\mathcal R}$, which does not belong to $S\cup T$, is
projected into some vector $r_1\in{\mathcal R}\cap S$. Therefore,
we have $r-r_1\in T$. On the other hand, $r-r_1$ belongs to the
flat $\mathbb R r+\mathbb R r_1$. Since ${\mathcal U}({\mathcal
R})$ is closed under intersection, the line $(\mathbb R r+\mathbb
R r_1)\cap T$ is an one dimensional flat of ${\mathcal
U}({\mathcal R})$, and, hence, there exists a vector
$r_2\in{\mathcal R}$ which spans this flat.

Now we consider the $3$-dimensional subspace $S+\mathbb R r_2$ of
$V$ and the unimodular system ${\mathcal R}\cap (S+\mathbb R r_2)$.
Obviously, the pure system of this unimodular system is closed under
intersection. Therefore, ${\mathcal R}\cap (S+\mathbb R r_2)$ is
isomorphic to $\mathbb A_2\oplus\mathbb A_1$. Thus, there can be at
most one generator outside of ${\mathcal R}\cap S$: the vector
$r_2$. However, we have another one: the vector $r\neq \pm r_2$. A
contradiction. Therefore ${\mathcal R}\subset S\cup T$ and the claim
is proven.\smallskip

Finally, suppose that ${\mathcal U}({\mathcal R})$ contains no flats
isomorphic to $\mathbb A_2$.  In such a case, we assert that
${\mathcal R}$ equals the sum of $n$ ($=\dim V$) exemplars $\mathbb
A_1$. Let $r_1,\ldots , r_n$ be linear independent elements of
${\mathcal R}$. We show that there holds ${\mathcal R}=\{\pm
r_1,\ldots , \pm r_n\}$. Assume some $r\in {\mathcal R}\setminus
\{\pm r_1,\ldots , \pm r_n\}$. Clearly we may assume that there
holds $r=r_1+\ldots + r_n$ (i.e. $r$ does not belong to the
coordinate hyperplanes). Let us consider the intersection of flats
$\mathbb R r_1+\mathbb R r_2$ and $\mathbb R r+\mathbb R r_3+\ldots
+\mathbb R r_n$. This intersection is a line $\mathbb R (r_1+r_2)$
and it is a flat of ${\mathcal U}({\mathcal R})$. Therefore, we have
$r_1+r_2\in{\mathcal R}$ and, hence, $\{\pm r_1,\pm r_2,\pm
(r_1+r_2)\}\subset {\mathcal R}$, but $\{\pm r_1,\pm r_2,\pm
(r_1+r_2)\}$ is isomorphic to $\mathbb A_2$. A contradiction.
$\Box$\medskip

Of course, the largest possible DC-classes are of the most
interest. Such DC-classes correspond to maximal pure systems and
maximal unimodular systems.\medskip

{\bf Definition.} A pure system ${\mathcal U}$ in $M$ is said to be
{\em maximal} if for any subspace $F$, not of ${\mathcal U}$, the
system ${\mathcal U}\cup \{F\}$ is not a pure system. A unimodular
system $\mathcal{R}$ is {\em maximal} if for any $r\notin
\mathcal{R}$ the system of vectors $\mathcal{R}\cup r$ is not a
unimodular.\medskip

{\bf Example 9.} {\em The unimodular system $\mathbb{A}_n$ is
maximal.} Let us remind a proof. Suppose that $r=(r_1,..., r_n)$
is an integer vector such that $\mathbb{A}_n \cup r$ is a
unimodular system. Since $\mathbb{A}_n$ contains the basic system
$\{\pm e_i, \ i=1,...,n\}$, all $r_i$ are equal to $0$ or $\pm 1$.
We assert that for any different $i$ and $j$ $r_ir_j=0$ or $-1$.
Indeed, suppose that $r_ir_j=1$. Let us consider the Abelian
subgroup $S$ generated by $r$, $e_i-e_j$, and all $e_k$, where
$k\ne i,j$. The index of $S$ in $\mathbb{Z}^n$ is equal to the
determinant of the matrix
$$\left(\begin{array}{cc}
  r_i & 1 \\
  r_j & -1
\end{array}\right),$$
that is $\pm 2$. This contradicts the purity of $S$. Therefore
only two of $r_i$ can differ of $0$ and in such a case these
coordinates are of opposite signs. That is $r\in \mathbb{A}_n$.

Let us reformulate this statement. Suppose that $L$ is a (rational)
line in $\mathbb{R}^n$ and $\rho$ is the canonical projection of
$\mathbb{R}^n$ onto $V'=\mathbb{R}^n/L$ such that the image
$\rho(\mathbb{A}_n)$ is a unimodular system in $V'$ (with respect to
the  integer structure $\rho(M)$). Then $L$ is generated by some
$r\in \mathbb{A}_n$ and the unimodular system $\rho(\mathbb{A}_n)$
is isomorphic to $\mathbb{A}_{n-1}$.

Indeed, unimodularity of $\rho(\mathbb{A}_n)$ means that $L$ is
mutually pure with any flat of $\mathbb{A}_n$. Hence $r$ belongs
to $\mathbb{A}_n$. The second assertion follows by considering of
the image of a subsystem $\mathbb{A}_{n-1}$ which is transversal
to $L$.\smallskip

We assert that $\mathbb{A}_n$ is not only maximal as a unimodular
system but also {\em the corresponding pure system
$\mathbb{U}=\mathcal{U}(\mathbb{A}_n)$ is maximal}. For this we
consider a vector subspace $F$ and suppose that $F$ is mutually pure
with any flat of $\mathbb{U}$. We have to show that $F$ also is a
flat of $\mathbb{U}$.

Let us consider the canonical projection $\phi$ of $\mathbb{R}^n$
onto the vector space $V'=\mathbb{R}^n/F$. As above, the image
$\phi(\mathbb{A}_n)$ is a unimodular system in $V'$. Let now $k$,
$1\le k \le n$, be a number such that $R=\mathbb{R}^k \cap F$ is
an one-dimensional subspace. Let us consider the restriction of
$\phi$ to $\mathbb{R}^k$. Since the image of $\mathbb{A}_k$ is a
unimodular set (as a subset of a unimodular set
$\phi(\mathbb{A}_n)$), we conclude that $R$ is generated by some
non-zero vector $r\in \mathbb{A}_k \subset \mathbb{A}_n$. Thus, we
have proven that $F$ contains some root $r$ of $\mathbb{A}_n$.

Now we consider the projection $\rho$ of $\mathbb{R}^n$ along
$\mathbb{R}r$. For the space $\mathbb{R}^n/\mathbb{R}r$ we have a
similar situation: a unimodular system
$\mathcal{R}=\rho(\mathbb{A}_n)$, isomorphic to
$\mathbb{A}_{n-1}$, and a vector subspace $F'=\rho(F)$ which is
mutually pure with flats of $\mathcal{R}$. By induction, $F'$ is a
flat of $\mathcal{R}$. Therefore its pre-image $F$ is a flat of
$\mathbb{A}_n$.\medskip

As a consequence, we obtain that the DC-class
$\mathcal{P}h(\mathbb{A} _n, \mathbb{Z} )$ of integer
g-poly\-matroids is maximal. $\Box$\medskip

{\bf Example 10.} The unimodular system $\mathbb{E}_5$ is maximal
too. However, the corresponding pure system is not maximal. In order
to see this, consider the following homomorphism $\phi: \mathbb{Z}^5
\to \mathbb{Z}$, $\phi(x_1,...,x_5)=x_1+...+x_5$. It is clear that
$\phi(r)=\pm 1$ for any root $r\in \mathbb{E}_5$. Therefore the
kernel of $\phi$, that is the hyperplane $H=[x_1+...+x_5=0]$, is
mutually pure with any flat of $\mathbb{R}_5$.

One can show that the S-class $\mathcal{P}h(\mathbb{E}_5)$ consists
of {\em zonohedra}, that is the sum of segments (bounded or not)
every of which is parallel to some root $r\in \mathbb{E}_5$. We
obtain that the intersection of two such integer zonohedra,
or the intersection of a
zonohedron and the hyperplane $H$, is an integer polyhedron.

The pure system corresponding the maximal unimodular system
$\mathbb{D}(K_{3,3})$ also is not maximal. It can be expanded by
adding some two-dimensional subspace. $\Box$\medskip

${\mathbb A}_n$ is a unique maximal unimodular system of dimension
$\le 3$. In dimension $4$, besides $\mathbb A_4$, there is another
maximal unimodular system $\mathbb D (K_{3,3})$. In dimension $5$,
there are $4$ non-isomorphic maximal unimodular systems; there are
$11$ in dimension $6$. For more details,  we refer to the article
\cite{DG}, which contains a complete description of maximal
unimodular systems.

Let ${\mathcal R}$ be a unimodular system. Elements $r$ of
${\mathcal R}$ can be identified with morphisms of $\mathbb A_1$ to
${\mathcal R}$. Conversely, morphisms of ${\mathcal R}$ to $\mathbb
A_1$ are called {\em co-roots}. In other words, a co-root is a
homomorphism of groups $\phi:M \rightarrow \mathbb Z$ such that
$|\phi(r)| \le 1$ for any root $r \in {\mathcal R}$. The set of
co-roots is denoted by ${\mathcal R}^*$.

A polyhedron is an ${\mathcal R}$-{\em polyhedron} if every of its
face is parallel to some ${\mathcal R}$-flat.  Denote by ${\mathcal
P}h({\mathcal R},\mathbb Z)$ the S-class of integer
$\mathcal{R}$-polyhedra. A pseudo-convex set $X$ in $M$ is said to
be ${\mathcal R}$-{\em convex set} if $\text{co}(X)$ is a ${\mathcal
R}$-polyhedron.

\section{Dual DC-classes associated to unimodular systems }

Besides the S-class of $\mathcal{R}$-polyhedra, we can associate
to a unimodular system $\mathcal{R}$ a dual I-class
integer $*\mathcal{R}$-polyhedra
(in the dual vector space $V^*$).\smallskip

Let $\mathcal{R}$ be a unimodular set in $M$, and let
$\mathcal{U}=\mathcal{U}(\mathcal{R})$ be the corresponding pure
S-system  in $V$. A polyhedron $P$ in $V^*$ is called
$*\mathcal{R}$-{\em convex} (or $*\mathcal{R}$-polyhedron) if it
belongs to $\mathcal{P}h(\mathcal{U}^{\bot})$, that is any face of
it is orthogonal to some $\mathcal{R}$-flat. In other words, a
$*\mathcal{R}$-polyhedron is given by a system of linear
inequalities (where $p$ is a linear functional on $V$)
$$
p(r)\le a(r), \text{ where } r\in \mathcal{R} \text{ and } a(r)\in
\mathbb{R}\cup \{+\infty\}.
$$
The inverse is also true. If all numbers $a(r)$ are integer, the
corresponding polyhedron is integer. Indeed, since the class
$\mathcal{U}^{\bot}$-polyhedra is I-class, we have to prove that
every hyperplane $H_r(a)=\{p\in V^*, \ p(r)=a\}$, where $r\in
\mathcal{R}$ and $a\in \mathbb{Z}$, contains an integer point. But
this is a consequence of primitiveness of $r$ in $M$. (This is a
kind of the Hoffman-Kruskal theorem \cite{H-K}.)

Thus, the set of all integer $*\mathcal{R}$-polyhedra is an I-class
of discrete convexity. For example, the class from Example 4 is the dual I-class
corresponding to the unimodular system $\mathbb{A}_n$.\medskip

In order to ``visualize'' integer $*\mathcal{R}$-polyhedra, it is
convenient to use the notion of a dicing \cite{Erdl}. A dicing is
the following regular polyhedral decomposition of $V^*$. Let us
consider the following countable (but locally finite) collection
of hyperplanes $H_r(a)=\{p\in V^*, \ p(r)=a\}$, where $r\in
\mathcal{R}$ and $a\in \mathbb{Z}$. These hyperplanes cut the
space $V^*$ on connected parts, the {\em regions} of the dicing.
Regions are bounded sets if $\mathcal{R}$ is of full dimension.
The closure of any region, as well as any its face, is called a
{\em chamber} of the dicing. The set $\mathcal{D}(\mathcal{R})$ of
the chambers form a polyhedral decomposition of $V^*$, that is the
chambers intersect by their faces and cover the whole space $V^*$.
If $\mathcal{R}$ is of full dimension, then the nodes of the
dicing (that is $0$-dimensional chambers) are integer points of
$V^*$, i.e., are elements of $M^*$.

Each chamber of the dicing $\mathcal{D}(\mathcal{R})$ is an
integer $*\mathcal{R}$-polyhedron. Conversely, any integer
$*\mathcal{R}$-polyhedron is a union of chamber of
$\mathcal{D}(\mathcal{R})$. Thus, an integer
$*\mathcal{R}$-polyhedron is nothing but a convex set composed of
chambers.\medskip

{\bf Example 11}. Let us consider the {\em dicing star} ${\bf
St}({\mathcal R})$. It is composed from those chambers of the
dicing ${\mathcal D}({\mathcal R})$, which contain the origin $0$.
In order to establish the convexity of ${\bf St}({\mathcal R})$,
we show that
$$
{\bf St}({\mathcal R})=\{p\in V^*, \ r(p)\le 1, \text{ where } r\in
\mathcal{R}\}.
$$
For the time being, we call ${\bf St}'$ the polyhedron appearing on
the right hand of the formula. Obviously any chamber which contains
$0$, belongs to ${\bf St}'$. Hence ${\bf St}({\mathcal R})\subset
{\bf St}'$.

Conversely, let $p \in {\bf St}'\setminus {\bf St}({\mathcal R})$.
Assume we move from $p$ to $0$ along the segment $[0,p]$. At some
time $t$, $0<t<1$, the point $t p$ will be on the boundary of ${\bf
St}({\mathcal R})$. Hence, there exists $r \in {\mathcal R}$ with
$r(tp)=1$. This implies that $r(p)=1/t>1$, a contradiction.

From this description of ${\bf St}({\mathcal R})$ we see that
integer points of ${\bf St}({\mathcal R})$ are the co-roots of
$\mathcal{R}$,
$$
{\bf St}({\mathcal R})(\mathbb{Z})=\mathcal{R}^*.
$$
Reversely, ${\bf St}({\mathcal R})=\text{co}(\mathcal{R}^*)$.
$\Box$\medskip

The dual pure system
$\mathcal{U}^{\bot}=\mathcal{U}(\mathcal{R})^{\bot}$ has the
following structure. It consists of the hyperplanes-mirrors
$H_r(0)=\mathbb{R}r^{\bot}$ and all possible intersections of the
mirrors. As well as a dicing, the mirrors cut the space $V^*$ onto
a finite number of cones (the {\em cameras}) which constitute a
fan $\Sigma(\mathcal{R})$ or $\mathcal{R}^{\bot}$. One-dimensional
flats of $\mathbb{U}$ are called {\em crossings} as well as their
primitive generators from $M^*$. (Of course, the crossings exist
only if the unimodular system $\mathcal{R}$ is of full dimension.)
As an element of $M^*$, a crossing is a surjective homomorphism of
Abelian groups $\xi:M \to \mathbb{Z}$ such that the kernel of
$\xi$ is a flat of $\mathcal{R}$. Let us denote
$\mathcal{R}^{\vee}$ the set of crossings in $M^*$.\medskip

{\bf Lemma 4.}  ${\mathcal R}^{\vee} \subset {\mathcal
R}^*$.\medskip

Proof. If $\mathcal{R}$ is not of full dimension, the set
$\mathcal{R}^{\vee}$ is empty. Therefore we can assume that
$\mathcal{R}$ is of full dimension. Let $\xi$ be a crossing, that is
a surjective homomorphism $M \to \mathbb{Z}$. Since the kernel
$\xi^{-1}(0)$ of $\xi$ is a flat, the image $\xi(\mathcal{R})$ is a
unimodular system in $\mathbb{Z}$, that is $\xi$ is a co-root.
$\Box$\medskip

{\bf Remark}. As Example 9 shows, for $\mathcal{R}=\mathbb{A}_n$
we have the equality ${\mathcal R}^{\vee} = {\mathcal R}^*$. For
other unimodular systems (such as $\mathbb{E}_5$) the crossings
constitutes a proper subset of ${\mathcal R}^*$.\medskip

In general case, the set ${\mathcal R}^{\vee}$ is not a unimodular
system in $M^*$; see Theorem 2. However, if we can find a
unimodular system $\mathcal{Q}$ in ${\mathcal R}^{\vee}$ (we call
such $\mathcal Q$ a {\em laminarization} of ${\mathcal R}$), this brings us an
advantage. Namely, in such a case we can construct $\mathcal{R}$-polyhedra simply as
$*\mathcal{Q}$-polyhedra. That is to define them by systems of linear inequalities
$$
\{v\in V, \ \xi(v)\le a(\xi), \ \xi \in \mathcal{Q}\}
$$
with arbitrary "right parts" $a(\xi)$. Of course, when $a(\xi)$
are integer, the corresponding polyhedron is integer too. Let us
give a more precise realization of this idea.\medskip

{\bf Example 12} (see also \cite{FT}). A family ${\mathcal T}$ of
subsets of a finite set $N$ is called {\em laminar} if for any $A,B
\in {\mathcal T}$, either $A \subset B$, or $B \subset A$,  or $A
\cap B=\emptyset$. Without loss of generality we can assume that any
singleton belongs to $\mathcal{T}$.

Let ${\mathcal T}$ be a laminar family. We assert that the set
${\mathcal Q}=\{ \pm {\bf 1}_T$, $T \in {\mathcal T}\}$ is a
unimodular set in the space $\mathbb{R}^N$. That is ${\mathcal Q}$
is indeed a laminarization of the system $\mathbb{A}_N$. Since the orthogonal
hyperplanes $({\bf 1}_T)^{\bot}$ are $\mathbb{A}_N$-flats, we have to check that any
intersections of such hyperplanes also are $\mathbb{A}_N$-flats.

Let us recall (see Example 5) that an $\mathbb{A}_N$-flat has the
form
$$
F(A_1,...,A_k):=\{x\in (\mathbb{R}^N)^*, \ x(A_j)=0 \text{ for }
j=1,...,k\},
$$
where $A_1,...,A_k$ are disjoint subsets of $N$. (The codimension
of $F(A_1,...,A_k)$ is equal to the number of non-empty $A_j$-s.)
In particular, the hyperplane $({\bf 1}_T)^{\bot}$ is $F(T)$. Let
us show that the intersection of hyperplanes
$F(T_1)$,...,$F(T_k)$, where $T_J \in \mathcal{T}$, has the form
$F(A_1,...,A_k)$. For this we write $A_j$ explicitly. Namely,
$A_j$ is equal to $T_j$ minus the union of those of $T_i$ which
are contained in $T_j$. Indeed, using the laminarity of
$\mathcal{T}$, we can assume that the $T_i$-s do not intersect.
Therefore vanishing $x(T_j)$-s is equivalent to vanishing
$x(A_j)$-s.

In particular, for a laminar family ${\mathcal T}$ in $N$, the
polyhedron defined by the inequalities
$$
a(S)\le x(S)\le b(S),\quad  S \in {\mathcal T},
$$
is an $\mathbb A_N$-polyhedron for any functions $a,\, b:{\mathcal
T} \rightarrow \mathbb R \cup \{\infty\}$, and is an integer
$\mathbb A_N$-polyhedron for integer-valued $a$ and $b$. $\Box$

\section{Exterior description of $\mathcal U$-polytopes}

In this section we characterize support functions of ${\mathcal
U}$-polyhedra, where $\mathcal{U}$ is a pure system. As we know
support functions of base polyhedra are closely related to
submodularity. Because of this, support functions of ${\mathcal
R}$-polyhedra give rise to a generalization of
submodularity.\medskip

Recall that the support function of a (non-empty) closed convex set
$A\subset V$ is the function  $\phi (A; \cdot ):V^* \rightarrow
\mathbb R \cup \{+ \infty \}$ on the dual space $V^*$ defined by the
following formula
\begin{equation}\label{supportn1}
\phi (A;p)=\sup_{x\in A} p(x),\quad p\in V^*.
\end{equation}
Let us work in a setting with compact sets in order to avoid messing
up with infinite values. In this setting the support function is
defined on whole the space $V^*$ and is homogeneous and convex.
Conversely, every homogeneous convex function $f$ on $V^*$ is the
support function of the {\it subdifferential} of $f$,
\begin{equation}\label{cc2}
\partial (f):=\{ x\in V\,|\,  x(p)\le f(p)\,\,\forall \,p\in V^*\}.
\end{equation}
The set $\partial (f)$ is non-empty, convex, and compact; and the
operations $\phi$ and $\partial $ are dual: $\partial(\phi (A))=A$
and $\phi (\partial f)=f$ (see, for example, \cite{R}).

Support functions of polytopes are characterized by a ``piece-wise
linearity'' property. It is convenient to use a notion of fan here.

A {\em fan} (in $V^*$) is a finite collection $\Sigma$ of polyhedral
cones possessing the following three properties: a) the cones
$\sigma \in \Sigma$ cover $V^*$; b) every face of any $\sigma \in
\Sigma$ is also in $\Sigma$; c) the intersection of two cones of
$\Sigma$ is a face of each of them. For example, in the previous
section we have defined the fan $\Sigma  ({\mathcal R})$.

A convex function $f$ on $V^*$ is {\em compatible} with a fan
$\Sigma$ if $f$ is linear on every cone $\sigma$ from $\Sigma$. In
this case, it is easy to show that $\partial (f)$ is a polytope.
More precisely, let $\sigma$ be a full-dimensional cone of the fan
$\Sigma$; denote by $v_{\sigma}$ a (unique) linear function on the
space $V^*$, which coincides with $f$ on the cone $\sigma$. Then
$v_{\sigma}$ (being considered as an element of $V$) is a vertex of
the polytope $\partial (f)$. And all vertices of the polytope are of
that form. In particular, a polytope $P$ is integer if and only if
its support function $\phi(P,\cdot)$ has integer values in integer
points. However, in this section, we shall not deal with the
integer-valuedness.

The support function of any polytope $P$ is compatible with the
following fan ${\mathcal N}(P)$. Given a point $x\in P$, the
following cone in the dual space $V^*$
$$
\mbox{Con}^*(P,x)= \{p \in V^*, p(x) \ge p(y) \ \forall y\in P\}
$$
is said to be the {\em cotangent} cone to $P$ at $x$. The collection
of all cotangent cones $\mbox{Con}^*(P,x)$, $x\in P$, forms the {\em
cotangent fan} (or the {\em normal fan}) ${\mathcal N}(P)$ of the
polytope $P$. For example, the cotangent fan of the zonotope
$\sum_{r\in {\mathcal R}}\mbox{co}(\{-r, r \})$ coincides  with the
arrangement fan $\Sigma ({\mathcal R})$. Cones of normal fan
${\mathcal N}(P)$ one-to-one correspond to faces of $P$. Moreover,
they are orthogonal one to other.

In particular, this gives the following\medskip

{\bf Proposition 7.} {\em Let $\mathcal{U}$ be a pure system in $V$, and let
$P\subset V$ be a convex polytope. The following assertion are equivalent:

a) $P$ is a $\mathcal{U}$-convex polytope;

b) the normal fan ${\mathcal N}(P)$ consists of
$\mathcal{U}^{\bot}$-cones.} $\Box$\medskip

When a pure system ${\mathcal U}$ is generated by a unimodular
system $\mathcal{R}$, we can say a bit more. In this case there is
the finest $*\mathcal{R}$-convex fan $\Sigma(\mathcal{R})$. And a
polytope $P$ is $\mathcal{R}$-convex if and only if its support
function is compatible with the fan $\Sigma(\mathcal{R})$.

One can give also the following characterization of
$\mathcal{R}$-polytopes.\medskip

{\bf Proposition 9.} {\em A polytope $P$ is $\mathcal{R}$-convex if
and only if there exists a polytope $P'$ such that $P+P'$ is an
$\mathcal{R}$-zonotope.}\medskip

Proof. It is clear that any edge of $P$ is parallel to some edge of
$P+P'$. Therefore $P$ is an $\mathcal{R}$-polytope. This prove the
"if" part of the statement.

Conversely, let $P$ be a $\mathcal{R}$-polytope. Then the
arrangement fan $\Sigma(\mathcal{R})$ is a refinement of the normal
fan $\mathcal{N}(P)$. Since the normal fan of an
$\mathcal{R}$-zonotope is $\Sigma (\mathcal{R})$, the assertion
follows from the following\medskip

{\bf Lemma 5 \cite{Grue}.} {\em For polytopes $P$ and $Q$ the
following assertions are equivalent:

     a) ${\mathcal N}(Q)$ is a refinement of ${\mathcal N}(P)$,

     b) there exists a polytope $P'$ such that $P+P'=kQ$, for some
$k\ge 0$.} $\Box$\medskip

Assume now that ${\mathcal R}$ is a full-dimensional unimodular
system, and that ${\mathcal R}^{\vee}$ is the set of crossings in
$M^*$. A function $f$, compatible with the fan  $\Sigma(\mathcal R)$
is uniquely determined by its restriction on ${\mathcal R}^{\vee}$,
that is by the family of real numbers $(f(\xi), \ \xi \in {\mathcal
R}^{\vee})$. However, the values $f(\xi)$, $\xi\in {\mathcal
R}^{\vee}$ are not arbitrary. Being the restriction of a convex
function, they must satisfy some kind of ``submodularity''
relations. These relations may be divided into two groups. The first
group of relations addresses the functions' linearity on each cone
of the fan. The second group of the relations yields convexity. Let
us formulate these relations more explicitly:\medskip

     I. Suppose that crossings $\xi_1,\ldots, \xi_m \in {\mathcal
R}^{\vee}$ belong to a cone $\sigma \in \Sigma({\mathcal R})$. Then
any linear relation $\sum_i \alpha_i \xi_i = 0$ should imply the
similar relation $\sum_i \alpha_i f(\xi_i)=0$.

     Of course, if the cone $\sigma$ is simplicial (as in the
case of $\mathbb A_n$), these relations disappear.\medskip

     II. Suppose that we have two adjacent (full-dimensional) cones
$\sigma$ and $\sigma'$ of the fan, separated by a wall $\tau$. Let
$\tau$ be spanned by the crossings $\xi_1,\ldots, \xi_m$, and let
$\xi, \xi'$ be crossings from $\sigma, \sigma'$ respectively, which
do not belong to the wall $\tau$. Then any relation $\alpha
\xi+\alpha' \xi'=\sum_i \alpha_i \xi_i$, where $\alpha, \alpha' >
0$, implies the relation $\alpha f(\xi)+\alpha' f(\xi') \ge \sum_i
\alpha_i f(\xi_i)$.

     According to Lemma 4, we can assume that $\alpha=
\alpha'=1$. But all the same, these relations do not look too
inspiring. In effect, it is neither easy to provide a collection of
numbers $(f(\xi), \xi \in {\mathcal R}^{\vee})$ satisfying the
relations I and II, nor easy to check that a given collection of
numbers satisfies these relations. See, nevertheless, a subsection
about laminarization.

Let us illustrate the above said for the unimodular systems
$\mathbb{A}(N)$ and $\mathbb{A}_N$.\medskip

{\bf Example 13. Base polytopes.} We show here that the class
${\mathcal B}$ of base polytopes (see Example 3) coincides with the
class of $\mathbb A(N)$-polytopes (a similar assertion is also true
for polyhedra; a proof, however, would involve support functions
with infinite values),  where $\mathbb A (N)$ is the unimodular
system from Example 6.

Recall that the set $\mathbb A(N)\subset (\mathbb R^N)^*$ consists
of differences $e_i-e_j$,  $i,j\in N$. Consider now how the
arrangement fan $\Sigma:=\Sigma(\mathbb A(N))$ in the space
$\mathbb R^N$ of functions on $N$ looks like. Given the root
$r=e_i-e_j$, the corresponding mirror $r^{\bot}$ consists of
functions $p \in \mathbb R^N$ satisfying the relation $p(i)=p(j)$.
This mirror divides the space of functions in two halfspaces
$\{p\,:\,p(i) \ge p(j)\}$ and $\{p:\,p(i) \le p(j)\}$. We see that
cones of the fan $\Sigma$ correspond to (weak) orders on $N$. If
$\preceq$ is an order, then the corresponding cone
$\sigma(\preceq)$ consists of monotone functions $p:(N,\preceq
)\rightarrow (\mathbb R,\le)$. For  example, full-dimensional
cones of $\Sigma$ correspond to linear orderings; the line of
constant functions $\mathbb R {\bf 1}_N$ corresponds to the total
indifference relation on $N$.

The set $\mathbb A(N)$ has full dimension in the hyperplane
$[x(N)=0]$ orthogonal to the constant function ${\bf 1}_N \in
\mathbb{R}^N$. Therefore we should consider the fan $\Sigma$ in the
factor space $\mathbb{R}^N/\mathbb{R}{\bf 1}_N$. The crossings
correspond  to dichotomous orders on $N$, which splits $N$ into two
classes $S$ and $N \setminus S$ ($S$ is different from $\emptyset$
and $N$). Therefore, crossings have the form ${\bf 1}_S$, $S\ne
\emptyset, N$.

Let now $f$ be a convex function compatible with the fan $\Sigma$.
Define the set-function $b:2^N \rightarrow \mathbb R$, $b(S)=f({\bf
1}_S)$ for $S \subset N$. We assert that $b$ is submodular. Indeed,
let $S$ and $T$ be subsets of $N$. Then, by convexity of $f$,
$$  b(S)+b(T)=f({\bf 1}_S)+f({\bf 1}_T)
\ge 2f(({\bf 1}_S+{\bf  1}_T)/2).   $$ On the other hand, since
$S\cap T \subset S\cup T$, the points ${\bf 1}_{S \cap T}$ and ${\bf
1}_{S \cup T}$ belong to a cone of $\Sigma$, and therefore
$$
b(S
\cap T)+b(S \cup T)=f({\bf 1}_{S \cap T})+f({\bf 1}_{S \cup T})=
2f(({\bf 1}_{S \cap T}+{\bf 1}_{S \cup T})/2).
$$
Since  $ {\bf 1}_S+{\bf 1}_T={\bf 1}_{S \cap T}+{\bf 1}_{S \cup T}$,
we have
$$
b(S)+b(T) \ge b(S \cap T)+b(S \cup T),
$$
that is $b$ is submodular function.

Conversely, any set-function $b$, considered as a function on the
set of vectors $\{{\bf 1}_S$, $S\subset N$\}, has the unique
extension $f={\tilde b}$ on whole $\mathbb R^{N}$ compatible with
the fan $\Sigma$. This extension coincides with the Choquet integral
(see \cite{Ch}) of the non-additive measure $b$, ${\tilde b}(p)=\int
pdb$. If $b$ is submodular function then ${\tilde b}$ is convex (see
\cite{L}).

The corresponding polytope $\partial {\tilde b}$ is given by the
following system of inequalities
$$
{\bf 1}_S(x)=x(S)\le b(S),\quad S\subset N,\quad x(N)=b(N), $$ and
is a base polytope. Thus, we prove\medskip

{\bf Proposition 10.} {\em  The class ${\mathcal P}t(\mathbb A(N))$
of $\mathbb A(N)$-polytopes coincides with the class of base
polytopes.}\medskip

Of course, the class of $\mathbb A(N)$-polyhedra coincides with the
class of base polyhedra, and the class of integer $\mathbb
A(N)$-polyhedra coincides with the class of integer base polyhedra.
$\Box$\medskip

{\bf Example 14. Generalized polymatroids.} In the same spirit, we
can check that the class of generalized polymatroids in $(\mathbb
R^N)^*$ coincides with the class of $\mathbb A_N$-polyhedra. The
arrangement ${\mathcal A}({\mathbb A_N})$ consists of hyperplanes
$p(i)=0$, $i\in N$, and $p(i)=p(j)$, $i,j \in N$. The collection of
vectors $\{\pm {\bf 1}_S$, $S\subset N\}$ is the set of crossings.
Cones of $\Sigma(\mathbb A_{n})$ are in a one-to-one correspondence
with pairs of orders $(\preceq_W , \preceq_{W'})$ on partitions
$(W,W')$ of $N$. These partitions derive from the partitions of
coordinates in non-negative and negative parts; $W$ denotes the
non-negative coordinates of vectors of a cone, whereas $W'$ denotes
the negative ones.

Now let $f$ be a convex function on $\mathbb{A}_N$ compatible with
the fan $(\Sigma (\mathbb A_N))$. Consider the following two
functions $a$ and $b$ on $2^N$: $a(S):=-f(-\mbox{\bf 1}_S)$ and
$b(S):=f(\mbox{\bf 1}_S)$ for $S \subset N$. There are three kinds
of relations between crossings: ${\bf 1}_S+{\bf 1}_T={\bf 1}_{S\cup
T}+{\bf 1}_{S\cap T},$ $-{\bf 1}_S-{\bf 1}_T=-{\bf 1}_{S\cup T}-{\bf
1}_{S\cap T},$ and
\begin{equation}\label{ident1} \mbox{\bf
1}_S +(-\mbox{\bf 1}_T)=\mbox{\bf 1}_{S-T}  +(-\mbox{\bf 1}_{T-S}).
\end{equation}
The first two yield submodularity of $b$ and supermodularity of
$a$, respectively, while the third yields the following
inequalities
\begin{equation}
b(S)-a(T)=f(\mbox{\bf 1}_S) +f(-\mbox{\bf 1}_T)\ge
f(\mbox{\bf 1}_{S-T})
+f(-\mbox{\bf 1}_{T-S})=b(S-T)-a(T-S).
\end{equation}
Thus, the pair $(b,a)$ is a strong pair in the sense of \cite{FT}. The
corresponding polyhedron $\partial f$ is given by the inequalities
$$
a(S) \le x(S) \le b(S),
$$
where $S \subset N$ and, by definition, $\partial f$ is a
generalized polymatroid.

Conversely, we can extend any strong pair $(b,a)$ to a convex
function on $\mathbb{R}^N$ compatible with the fan $\Sigma(\mathbb
A_N)$. Thus, the class of (bounded) generalized polymatroids
coincides with the class of $\mathbb{A}_N$-polytopes. Similarly, the
class of all generalized polymatroids coincides with the class of
$\mathbb{A}_N$-polyhedra, and the class of integer generalized
polymatroids coincides with the class of integer
$\mathbb{A}_N$-polytopes.

\section*{Acknowledgements}
We thank Christine Lang and anonymous referees for helpful comments and remarks.

\noindent 
Central Institute of Economics and Mathematics\\ Nahimovski Prospect 47,
Moscow 117418, Russia\\ e-mail: vdanilov43@mail.ru\\ koshevoy@cemi.rssi.ru

 \end{document}